\documentclass[11pt]{amsart}
\usepackage{geometry}                
\usepackage{graphicx}
\usepackage{amssymb}
\usepackage{epstopdf}
\newtheorem{defn}{Definition}
\newtheorem{thm}{Theorem}

\newcommand{\eps}{\varepsilon}
\newcommand{\R}{\mathbb{R}}
\newcommand{\calK}{\mathcal{K}}
\newcommand{\calM}{\mathcal{M}}

\title{Convergence to equilibrium 
for a class of exchange economies}
\author[R.S.MacKay]{R.S.MacKay\\
Mathematics Institute \& Centre for Complexity Science, \\University of Warwick, Coventry CV4 7AL, U.K.\\
R.S.MacKay@warwick.ac.uk}

\date{\today}                                           

\begin{document}

\begin{abstract}
For a class of stochastic dynamical models of exchange economies that we call ``fully connected Cobb-Douglas'', the paper proves convergence of the probability distribution to an equilibrium, in total variation metric as time goes to infinity.  The convergence is exponential and the equilibrium is determined uniquely by the number of agents, their ``exponents'', and the initial amounts of money and goods in the economy.
\end{abstract}

\footnote{ORCID:0000-0003-4771-3692} 

\maketitle

\section{Introduction}
For this paper, an exchange economy consists of a number $N\ge 2$ of agents who exchange real amounts of a finite number $M$ of types of good in pairwise encounters according to Markovian dynamics described in the following three paragraphs.  

Pairs of agents $i,j$ make encounters independently at rate $k_{ij}\ge 0$ for some symmetric $N\times N$ matrix $k$ with $k_{ii}=0$.  For this paper, we take the encounter matrix to be ``fully connected'', meaning $k_{ij}>0$ for all $i\ne j$.


Each agent $i$ has a ``utility function'' $u_i$, which is a non-negative function of the vector $g$ of (non-negative) amounts $g_1,\ldots g_M$ of each type of good that it  owns.  For this paper we will take what we call ``Cobb-Douglas'' agents, for which $u_i$ has the form
\begin{equation}
u_i(g) = \prod_{m=1}^M g_m^{\alpha_m -1}
\end{equation}
for some ``exponents'' $\alpha_m>0$, that may also depend on $i$.

At each encounter between agents $i,j$, they pool their belongings to make a vector $g \in \R_+^M$ ($\R_+ = \{x \in \R:~x\ge 0\}$) of amounts of the types of good, and then redistribute $g$ between them to make as outcome a pair of vectors $(g^i, g^j)$ (we use superscripts here to indicate agents, as opposed to powers, to avoid confusion with subscripts for type of goods) with probability density proportional to
\begin{equation}
u_i(g^i)u_j(g^j)
\end{equation}
on the affine space defined by $g_m^i \ge 0$, $g_m^j \ge 0$, $g_m^i + g_m^j = g_m$ for all $1 \le m \le M$.
The redistribution is chosen independently of all previous events.

For a number $x \in \R^+ = \{ x \in \R:~x>0\}$, the {\em $(N-1)$-simplex} $\Delta_{N,x}$ is the set of vectors $(g_1,\ldots g_N) \in \R_+^N$ with $g_i\ge 0$ and $\sum_{i=1}^N g_i = x$ (here we have switched to using subscripts for agents, as we think of $x$ as representing the total amount of a single type of good).
Given the vector $G \in \R_+^M$  of initial total amounts of goods, the result of the above three paragraphs is a Markov process on the product of simplices $\prod_{m=1}^M \Delta_{N,G_m}$.

The {\em Dirichlet distribution} $D(\alpha,G)$ (for a vector $\alpha \in {\R^+}^N$ and a single type of good) is the probability distribution on $\Delta_{N,G}$ with density (labelling agents by subscripts now)
\begin{equation}
\frac{1}{Z(\alpha,G)} \prod_{i=1}^N g_i^{\alpha_i -1} 
\label{eq:Dirichlet}
\end{equation}
and normalisation constant
\begin{equation}
Z(\alpha,G) =  \frac{\prod_{i=1}^N \Gamma(\alpha_i)}{\Gamma(s_N)} G^{s_N},
\end{equation}
where $\Gamma$ is the Gamma-function and
\begin{equation}
s_N = \sum_{i=1}^N \alpha_i.
\label{eq:sN}
\end{equation}

The process is reversible with respect to the probability distribution
\begin{equation}
\prod_{m=1}^M D(\alpha_m,G_m),
\end{equation}
where $G_m$ represents the vector of amounts of good of type $m$ for the agents and $D(\alpha_m,G_m)$ is the Dirichlet distribution on the simplex $\Delta_{N,G_m}$ for exponents $\alpha_m$ (which is in general a vector of values indexed by the agents).  This can be seen from symmetry of the probability flux between any two states when in this distribution.  It follows that the distribution is an equilibrium.  Note that it is independent of the matrix $k$.

The goal of this paper is to prove 

\begin{thm}
\label{thm:main}
For any initial probability distribution on the product of simplices, the probability distribution for the process at time $t>0$ converges to the product of Dirichlet distributions in total variation metric as $t \to \infty$.
\end{thm}

\noindent The definition of total variation metric is recalled in the next section.

The strategy of proof is described in Section~\ref{sec:strategy}.  The proof of the key step of the strategy is given in Section~\ref{sec:full}.  
The paper concludes with a discussion in Section~\ref{sec:disc}.

There are some precursors in the literature.  Carlen et al \cite{CCL} study the Kac model for velocities $v_i$ in a 1D gas, which is equivalent to a single type of good (kinetic energy $E_i = \tfrac12 v_i^2$) with exponent $\alpha = \frac12$ (uniform density on the circle $v_i^2 + v_j^2 = 2$ maps to Beta$(\tfrac12,\tfrac12)$ on the simplex $E_i+E_j=1$) and identical encounter rates between all pairs.  They prove exponential convergence to equilibrium in a Hilbert space and determine the spectrum of decay rates.
D\"uring et al \cite{D+} study the case of a single type of good with $\alpha=1$ and identical encounter rates, and obtain convergence to equilibrium in total variation metric but by a procedure that we consider over-complicated and not completely clear.  Nonetheless, their approach by induction on $N$ is valuable and we imitate that here.

\section{Strategy}
\label{sec:strategy}

The process treats each type of goods independently, so it is enough to consider a single type of good.
 We suppose $N\ge 2$ agents and total amount $G>0$ of the good.

The times of encounters form a Poisson process with rate $K = \sum_{i<j} k_{ij}$, i.e.~the ``waiting'' time between encounters is i.i.d.~$Exp(K)$.  So with probability 1, the encounter times are distinct and have no accumulation point.  For $t>0$ the probability that the number $m$ of encounters in $[0,t]$ is 
\begin{equation}
p_m(t) = e^{-Kt}\frac{(Kt)^m}{m!}.
\end{equation}

Let $\tilde{P}$ be the transition operator on $\Delta_{N,G}$ that selects a pair $(i,j)$ of agents with probability $\frac{k_{ij}}{K}$ and redistributes the goods between them with probability density proportional to $u_i(g_i)u_j(g_j)$.  The kernel $\tilde{P}_x$ for initial vector $x$ of amounts of good indexed by all agents is a probability distribution supported on $\left(\begin{array}{c} N \\ 2 \end{array}\right)$ lines through $x$, parallel to the edges of $\Delta_{N,G}$, with weights $k_{ij}/K$ and densities proportional to $u_i u_j$.

Then the transition operator for time $t>0$ for the full process is
\begin{equation}
P^t = \sum_{m\ge 0} e^{-Kt} \frac{(Kt)^m}{m!} \tilde{P}^m.
\end{equation}

Recall the definition of total variation metric.

\begin{defn}
The {\em total variation distance} between two probability measures $\mu,\nu$ on a measure space $\calM$ is
\begin{equation}
d_{TV}(\mu,\nu) = \sup_{A} \mu(A)-\nu(A)
\end{equation}
over measurable subsets $A$ of $\calM$.
\end{defn}

Say two non-negative measures $\mu, \nu$ on a measure space $\calM$ satisfy $\mu \ge \nu$ if for all measurable subsets $A$ of $\calM$, $\mu(A) \ge \nu(A)$.

The strategy of proof of Theorem~\ref{thm:main} is to apply Doeblin's theorem:

\begin{thm}[Doeblin]
If $P$ is a transition operator on a measure space $\calM$ and there is a non-negative measure $\beta$ on $\calM$ with $\beta(\calM)= \eps>0$ such that for all $x \in \calM$, $P_x \ge \beta$, then for all probability measures $\mu, \nu$ on $\calM$,
\begin{equation}
d_{TV}(\mu P, \nu P) \le (1-\eps) d_{TV}(\mu,\nu).
\end{equation}
In particular, $P$ has a unique fixed point $\pi$ and it is exponentially attracting in total variation metric: $d_{TV}(\mu P^n,\pi) \le (1-\eps)^n d_{TV}(\mu,\pi)$ for $n\ge 0$.
\end{thm}

Some call this version, which is for a general measure space, Harris' theorem, e.g.~\cite{BH}.  For a proof, see the Appendix.

We will construct a non-negative measure $\beta_N$ on $\Delta_{N,G}$ with $\beta_N(\Delta_{N,G})>0$ such that for all initial states $x$, $\tilde{P}_x^{N-1} \ge \beta_N$.  This will be done by induction on $N$.  It follows that for any $t>0$, 
$$P_x^t \ge e^{-Kt} \frac{(Kt)^{N-1}}{(N-1)!} \beta_N .$$  
Choose a $\tau>0$. Then by Doeblin's theorem, $P^\tau$ has a unique fixed point and it is exponentially attracting.  We already know by reversibility that the Dirichlet distribution is a fixed point, call it $\pi$.  So $d_{TV}(\mu P^{n\tau},\pi) \to 0$ exponentially as $n \to \infty$.  It is standard that for any transition operator, $d_{TV}(\mu P^t, \nu P^t)$ is nowhere increasing in $t$, so we can fill in the times between the $n\tau$.  So Theorem~\ref{thm:main} follows.

Actually, we will obtain $\tilde{P}_x^{N-1} \ge \beta_N = c_N D(\alpha,G)$ for some constant $c_N\le 1$, and $D(\alpha,G)$ is invariant under $\tilde{P}$ so $\tilde{P}_x^m \ge \beta_N$ for all $m \ge N-1$, and hence $$P^\tau \ge e^{-K\tau} \sum_{m\ge N-1} \tfrac{(K\tau)^m}{m!} \ c_N D(\alpha,G),$$ which strengthens the estimate.  In particular, $$\eps = c_N e^{-K\tau} \sum_{m\ge N-1} \tfrac{(K\tau)^m}{m!}.$$  The resulting convergence rate with respect to time (as opposed to iteration of $\tilde{P}$) is at least $\tau^{-1} \log \frac{1}{1-\eps}$.  This can be maximised over the choice of $\tau$ if desired, but the true convergence rate is likely to be much larger.

\section{Construction of $\beta_N$}
\label{sec:full}

The matrix $k$ is finite so there is $\calK>0$ such that for all $i \ne j$, $k_{ij} \le \calK$.
Because $k$ is fully connected  there is $\kappa>0$ such that $k_{ij} \ge \kappa$ for all $i\ne j$.  Let $\rho = \kappa/\calK \in (0,1]$.

We will show that for all states $x \in \Delta_{N,G}$, 
\begin{equation}
\tilde{P}_x^{N-1} \ge \beta_N = c_N D(\alpha,G)
\label{eq:ind}
\end{equation}
for some $c_N>0$ to be determined (depending on $\rho$ and the vector $\alpha \in {\R^+}^N$).

For $N=2$ we can take $c_2 = 1$, because $\tilde{P}_x = D(\alpha_1,\alpha_2,G)$ for all $x = (g_1,g_2)$ with $g_1+g_2=G$.  

Suppose the hypothesis (\ref{eq:ind}) is satisfied for some number $N\ge 2$ agents.

Take $N+1$ agents.  Because the total amount of goods is $G$, at least one agent has initial goods $y \le \frac{G}{N+1}$.  Label such an agent by $N+1$.  The remaining amount of goods $\tilde{G} = G -y \ge \frac{N}{N+1}G$, and $(g_1,\ldots g_N) \in \Delta_{N,\tilde{G}}$.
We obtain a lower bound for $\tilde{P}_x^{N-1}$ by ignoring encounters with agent $N+1$ and using the induction hypothesis (\ref{eq:ind}).
The share of probability in $\tilde{P}_x$ that goes to encounters with agent $N+1$ at each iteration is at most 
$$\frac{N\calK}{N\calK + \tfrac12 N(N-1)\kappa} = (1+ \tfrac12(N-1)\rho)^{-1}.$$
So the complementary probability is at least
$$1-(1+\tfrac12 (N-1)\rho)^{-1} = (1+ \tfrac{2}{(N-1)\rho})^{-1}.$$
Thus
\begin{equation}
\tilde{P}_x^{N-1} \ge \left(1+ \tfrac{2}{(N-1)\rho}\right)^{1-N} c_N D(\alpha,\tilde{G})
\end{equation}
on $\Delta_{N,\tilde{G}}$.

Then we apply $\tilde{P}$ once more, this time bounding it from below by taking only encounters with agent $N+1$.  The share of probability for each of these is at least $\tfrac{2\rho}{N(N+1)}$.
We obtain
\begin{equation}
\tilde{P}_x^N \ge \left(1+ \tfrac{2}{(N-1)\rho}\right)^{1-N}  \tfrac{2\rho}{N(N+1)} \, c_N D(\alpha,\tilde{G}) \otimes \delta_{y} \sum_{i\le N} E_i
\label{eq:12}
\end{equation}
on $\Delta_{N+1,G}$,
where $\delta_y$ represents the atomic distribution on $g_{N+1}=y$ and $E_i$ is the transition operator for the encounter of $i$ with $N+1$.
It is enough to lower bound this in the subset of $\Delta_{N+1,G}$ where $g_1 \ge g_j$ for $j = 2 \ldots N$, the other similar subsets being treated the same. For this we can use just the term $i=1$ in the sum.  

The kernel of the transition operator $E_1$ from $g$ to $g'$ has ``density'' (using delta-functions to keep track of the restrictions to affine subspaces) 
\begin{equation}
\frac{{g_1'}^{\alpha_1-1} {g_{N+1}'}^{\alpha_{N+1}-1}} {Z(\alpha_1,\alpha_{N+1},g_1+g_{N+1})}   \delta(g_1'+g_{N+1}' - g_1 - g_{N+1}) \prod_{i=2}^N \delta(g'_i - g_i).
\end{equation}
So the product distribution $D(\alpha,\tilde{G}) \otimes \delta_{y}\, E_1$ has density as a function of $g'$
\begin{align}
\label{eq:integral}
\int & \frac{ \prod_{i=1}^N  g_i^{\alpha_i-1} }{Z(\alpha_1,\ldots \alpha_N, G-y)}\, \delta\left(\sum_{i=1}^N g_i - G + y\right) \delta(g_{N+1}-y)  \times \\
& \frac{{g_1'}^{\alpha_1-1} {g_{N+1}'}^{\alpha_{N+1}-1} }{Z(\alpha_1,\alpha_{N+1},g_1+g_{N+1})}\, \delta(g_1'+g_{N+1}' - g_1 - g_{N+1}) \prod_{i=2}^N \delta(g'_i - g_i) \prod_{i=1}^{N+1} dg_i. \nonumber
\end{align}
The integral over $g_{N+1}$ is achieved by substituting $g_{N+1}$ by $y$.
The integrals over $g_i$ for $i=2,\ldots N$ are achieved by substituting $g_i$ by $g'_i$.
The integral over $g_1$ is achieved by substituting $g_1$ by $g'_1+g'_{N+1} - y$.
Writing
\begin{equation}
x = g'_1 + g'_{N+1},
\end{equation}
we obtain
\begin{equation}
\frac{ (x-y)^{\alpha_1-1} \prod_{i=2}^N {g'_i}^{\alpha_i-1}} {Z(\alpha_1,\ldots \alpha_N,G-y)} \frac{{g'_1}^{\alpha_1-1} {g'_{N+1}}^{\alpha_{N+1}-1}}{Z(\alpha_1,\alpha_{N+1},x)}\, \delta\left(\sum_i g'_i - G\right),
\label{eq:int2}
\end{equation}
restricted to $x\ge y$.
Recognising a large part of this as the density for the Dirichlet distribution on $\Delta_{N+1,G}$ from (\ref{eq:Dirichlet}), it boils down to
\begin{equation}
\frac{(x-y)^{\alpha_{1}-1} }{Z(\alpha_1,\alpha_{N+1},x)} \frac{Z(\alpha_1,\ldots \alpha_{N+1}, G)}{Z(\alpha_1,\ldots \alpha_{N},G-y)} D_{N+1,G}(g')\, \delta\left(\sum_{i=1}^{N+1} g'_i - G\right)
\label{eq:density}
\end{equation}
for $x\ge y$.
Now 
\begin{equation}
Z(\alpha_1,\alpha_{N+1},x) = x^{\alpha_1+\alpha_{N+1}-1} \frac{\Gamma(\alpha_1)\Gamma(\alpha_{N+1})}{\Gamma(\alpha_1+\alpha_{N+1})}.
\end{equation}
Also
\begin{equation}
\frac{Z(\alpha_1,\ldots \alpha_{N+1}, G)}{Z(\alpha_1,\ldots \alpha_{N},G-y)} = \Gamma(\alpha_{N+1}) \frac{\Gamma(s_N)}{\Gamma(s_{N+1})} \frac{G^{\alpha_{N+1}}} {(1-\frac{y}{G})^{s_N} },
\end{equation}
using the notation from (\ref{eq:sN}).  Note that $\frac{N}{N+1} \le 1-\frac{y}{G} \le 1$ (of which we shall use only the upper bound).

Evaluate (\ref{eq:density}) in the subset where $g'_1 \ge g'_i$ for $i=2,\ldots N$.  Then $g'_1 \ge \frac{G-g'_{N+1}}{N}$.  So $x \ge \frac{G}{N}$.  We have also that $y \le \frac{G}{N+1}$, hence $x-y \ge \frac{G}{N(N+1)}$.  The function $h$ defined by 
\begin{equation}
h(x) = \frac{(x-y)^{\alpha_{1}-1} }{x^{\alpha_1+\alpha_{N+1}-1}}\, G^{\alpha_{N+1}}
\label{eq:h}
\end{equation}
is continuous and positive for $x$ in the interval $[y+\frac{G}{N(N+1)},G]$, so has a positive lower bound.  By scaling $x,y$ to $G$ we see that it is independent of $G$.  A positive lower bound can be chosen uniformly in $y \le \tfrac{G}{N+1}$.  Also, we can take the minimum over change of $\alpha_1$ to $\alpha_i$ for $i=2\ldots N$.   We denote the resulting lower bound by $L_N$.

Also, let $J_N$ be a lower bound for  
\begin{equation}
\frac{\Gamma(\alpha_1+\alpha_{N+1})}{\Gamma(\alpha_1)} \frac{\Gamma(s_N)}{\Gamma(s_{N+1})} 
\label{eq:J}
\end{equation}
over permutations of the agents.

Let 
\begin{equation}
c_{N+1} =  \left(1+ \tfrac{2}{(N-1)\rho}\right)^{1-N}  \tfrac{2\rho}{N(N+1)} J_N  L_N  \, c_N
\end{equation}
Then
\begin{equation}
\tilde{P}_x^N \ge c_{N+1} D(\alpha,G),
\end{equation}
so (\ref{eq:ind}) holds for $N+1$.  

By induction on $N$, this completes the proof of (\ref{eq:ind}) for all $N\ge 2$.  By the strategy of Section~\ref{sec:strategy}, this completes the proof of Theorem~\ref{thm:main}.

\section{Discussion}
\label{sec:disc}

The result of this paper simplifies and generalises that of \cite{D+}.  Furthermore, it
justifies Axiom A0 of a thermodynamic approach to macroeconomics \cite{CM} in the case of fully connected Cobb-Douglas economies.

It would be good to generalise the result to general connected Cobb-Douglas economies, in the sense that for all $i\ne j$ there is an $m>0$ and a finite path $i_0=i,\ldots , i_{m} = j$ with $k_{i_n, i_{n+1}} >0$ for $n = 0,\ldots m - 1$.  It was not obvious to us how to extend the use of Doeblin's theorem to this case.

It would be good also to extend the analysis to other utility functions, in particular not of the form of a product over types of good, and to allow some dependence of each agent's utility function on the amounts of other agents' possessions.  Further extensions would be to treat cases where agents make available only a fraction of their possessions at each encounter and other variants.

One might hope for an estimate of convergence rate uniform in number of agents (under suitable uniform hypotheses), but this is unlikely.  For example, for agents arranged in a one-dimensional nearest-neighbour chain, equilibration is analogous to diffusion, which takes a time proportional to the square of the length.  It might also be necessary to replace total variation metric (compare \cite{M} for probabilistic cellular automata).

Finally, it would be good to derive an order of magnitude estimate of the convergence rate (as opposed to the probably pessimistic rigorous bound we have obtained here).

\section*{Acknowledgements}
I am grateful to Bertram D\"uring, Enrico Scalas and Nicos Georgiu for answering questions about their paper \cite{D+}, and to Jack Lewis for going through my proof.

\appendix

\section{Proof of Doeblin's theorem}

Applying $P_x \ge \beta$ to $\calM$ implies $1 \ge \beta(\calM)$.
Similarly, for any probability measure $\mu$, $\mu P \ge \beta$.  

If $\beta(\calM)=1$ then it follows that $\mu P = \beta$.  Then for any other probability measure $\nu$, $\nu P = \beta$ too, so $d_{TV}(\mu P,\nu P) = 0$ and the result follows.

If $0< \eps =\beta(\calM)<1$ then for each $x \in \calM$ let
\begin{equation}
K_x = \frac{P_x-\beta}{1-\eps}.
\end{equation}
It is a Markov kernel (maps probabilities to probabilities).  Then for any probability measure $\mu$,
\begin{equation}
\mu P = \beta + (1-\eps)\mu K,
\end{equation}
so for any pair $\mu, \nu$,
\begin{equation}
d_{TV}(\mu P, \nu P) =(1-\eps) d_{TV}(\mu K, \nu K) \le (1-\eps) d_{TV}(\mu,\nu),
\end{equation}
where the inequality comes from the equivalent definition $d_{TV}(\mu,\nu) = \tfrac12 \sup_f (\mu f - \nu f)$ over measurable functions $f$ with  $|f| \le 1$.
This concludes the proof.

\end{document}